\documentclass[11pt]{amsart}

\usepackage{amssymb,amsmath,amsfonts,amsthm,amscd,eufrak}

\usepackage{graphics}
\usepackage[colorlinks=true]{hyperref}

\input xy
\xyoption{all}

\newcommand{\CC}{\mathbb C}
\newcommand{\RR}{\mathbb R}
\newcommand{\maA}{\mathcal{A}}

\newcommand{\maB}{\mathcal{B}}

\newcommand{\maD}{\mathcal{D}}

\newcommand{\gp}{\Gamma}

\newcommand{\ZZ}{\mathbb{Z}}

\newcommand{\maU}{\mathcal{U}}

\newcommand{\maV}{\mathcal{V}}

\newcommand{\pdoM}{\Psi^{\infty}(M)/\Psi^{-\infty}(M)}
\newcommand{\smooth}{\mathcal{C}^{\infty}}

\newcommand{\ip}[1]{\langle #1 \rangle}

\theoremstyle{definition}
\newtheorem{definition}{Definition}[section]

\theoremstyle{plain}
\newtheorem{theo}[definition]{Theorem}
\newtheorem{prop}[definition]{Proposition}
\newtheorem{lem}[definition]{Lemma}
\newtheorem{cor}[definition]{Corollary}

\theoremstyle{remark}
\newtheorem{remark}[definition]{Remark}

\begin{document}

\textcolor{red}{
 \title[Traces] {An equivariant non-commutative residue } 
}

\author[S Dave]{Shantanu Dave}
\address{University of Vienna, Austria }
\email{shantanu.dave@univie.ac.at}
\thanks{Supported by FWF grant Y237-N13 of the Austrian Science Fund.}

\begin{abstract}
Let $\gp$ be a finite group acting on a compact manifold $M$ and
$\maA(M)$ denote the algebra of classical complete symbols on $M$. We
determine all traces on the cross product algebra $\maA(M) \rtimes
\Gamma$.  These traces appear as residues of certain meromorphic
'zeta' functions and can be considered as equivariant generalization
of the non-commutative residue trace.  The local formula for these
traces depends on more than one component of the complete asymptotic
expansion. For instance, 
the local formula for these traces
depends also on derivatives in the normal directions to fixed point
manifolds of higher order components. As an application, we obtain a
formula for the asymptotic occurrence of an irreducible representation
of $\gp$ in the eigenspaces of an invariant positive elliptic
operator. We also obtain an new construction for Dixmier trace of an
invariant operator.
\end{abstract}

\maketitle

\section{Introduction}\label{sec_intro}
The spectrum of a geometric differential operator, for instance that
of the Laplace operator, is of significant interest in analysis and
geometry. It is well known that the Laplace operator on a bounded
domain has a discrete spectrum. In a classical paper \cite{Weyl}
Hermann Weyl derived estimates on asymptotic growth of eigenvalues. In
1980s, {Guillemin} \cite{Guillemin} has obtained Weyl's estimates by
using the non-commutative residue.

In \cite{Wodzicki},  Wodzicki introduced a trace on the algebra
of pseudodifferential operators on a smooth, compact manifold,
called the {\em residue trace}. A similar result was obtained in
\cite{Guillemin}.
Many interesting invariants can be expressed in terms of traces on
algebras. The space of all traces on an algebra forms its 0-th
Hochschild cohomology.  Hochschild and cyclic homology are important
tools in non-commutative geometry.  Results on the cyclic homology for
algebras of complete symbols over compact manifolds were obtained in
\cite{Br-Gz,Wodzicki}. In particular, these homology calculations
recover the non-commutative residue of {Guillemin} \cite{Guillemin}
and Wodzicki \cite{Wodzicki}. The Atiyah-Singer index theorem can also
be formulated in terms of the non-commutative residue. 
Residues form the basic ingredient in the non-commutative index
theorem of Connes and Moscovici \cite{Co-Mo}.
An homological approach to the equivariant version of the
non-commutative residue is discussed in \cite{Dave}

We construct an equivariant, delocalized generalization of the
noncommutative residue for manifolds endowed with a group action. This
result amounts, in particular, to the computation of the traces on the
cross-product algebra of pseudodifferential operators by a finite
group. As an application, we then obtain results on the asymptotic
properties of eigenspaces of operators that are invariant with respect
to the given group action. This application is to count the occurrence
of an irreducible representation asymptotically in the eigenspaces.

Let $V$ be a smooth, compact, $n$ dimensional Riemannian manifold, and
let $\Delta := d^*d$ be the (positive) Laplace operator on $V$. Weyl's
theorem asserts that the spectral counting function
\begin{equation*}
        N_{\Delta}(\lambda) := \#\{\lambda_i: \lambda_i\in
	\mathrm{sp}(\Delta), \lambda_i\leq\lambda\}
\end{equation*}
grows asymptotically as $C\lambda^{\frac{n}{2}}$, where the constant
$C$ depends only on the dimension and volume of $M$.

In view of the results of Guillemin and Wodzicki, in order to
generalize Weyl's result to the equivariant case, we can proceed as
follows.  Let $\Gamma$ be a finite group, acting on a closed manifold
$M$ by diffeomorphisms. For instance, $\Gamma$ could be a finite group
of isometries. Let $D$ be a positive differential operator on $M$ that
is preserved under the $\Gamma$ action. A good choice is the Laplace
operator $D=\Delta$ associated to an invariant metric. More generally,
$D$ could be an elliptic pseudodifferential operator of positive
order. Let $\{\lambda_i\}$ be the eigenvalues of $D$ and let
$V_i=\textrm{ker}(D-\lambda_iI)$ be the corresponding
eigenspaces. Then each eigenspace $V_i$ acquires a representation of
the group $\Gamma$. One natural problem is to asymptotically count the
occurrence of any particular representation in these eigenspaces. More
precisely, if we start with a representation $\pi$ of $\Gamma$ and we
count the multiplicity of $\pi$ in each eigenspace $V_i$, we obtain
the refined counting function
\[  N_{\pi,D}(\lambda):=\sum_{\lambda_i<\lambda}\
    \textrm{``multiplicity of}\ \pi\ \textrm{in}\ V_i.\textrm{''}
\]
Then a natural question is to estimate the asymptotic growth of
$N_{\pi,D}(\lambda)$.  One might also want to compare asymptotically
the relative occurrence of various representations of $\Gamma$ in
these eigenspaces. The answer that we propose to these questions is
based on the construction of ``equivariant, delocalized''
noncommutative residues.

Recall that the classical pseudodifferential operators on a closed
manifold $M$ form an algebra $\Psi^{\infty}(M)$. The space of
smoothing operators $\Psi^{-\infty}(M)$ is an ideal of this
algebra. The algebra of complete symbols is then the quotient
$\Psi^{\infty}(M)/\Psi^{-\infty}(M)$, and we will denote it by
$\maA(M)$.  If the cosphere bundle $S^*M$ is connected, then the
noncommutative residue is the unique (up to a constant) trace on the
algebra $\maA(M)$.  Moreover, if $P$ is an order one, positive
operator, then, for any $A\in \Psi^{\infty}(M)$, we can define
$\zeta_A(z):=Tr(P^{-z}A)$ for $\Re(z)$ large enough. This function
turns out to be meromorphic and the residue at zero,
$Tr_R(A)=\mathrm{res}_{z=0} \,\zeta_A(z)$, then defines a trace on
$\Psi^{\infty}(M)$, called the non-commutative residue
\cite{Guillemin, Wodzicki}. This trace descends to a trace on
$\maA(M)$.

To obtain an equivariant version, one should consider traces on the
cross-product algebra $\maA(M) \rtimes \Gamma$ obtained from an action
of $\Gamma$ on $M$. As a set, the cross product algebra is 
\begin{equation}
\maA(M) \rtimes \Gamma := \big\{\ \sum_{g\in\Gamma}A_g g :\ A_g\in
\maA(M)\ \big\}.
\end{equation}
Its product structure comes from the action of $\Gamma$, namely,
\begin{equation*}
A_g g \cdot B_h h = A_g g(B_h) gh.
\end{equation*}
For a fixed, $\Gamma$ invariant operator $D$ of positive order and any
element $g$ in $\Gamma$, we consider the function $z\mapsto
\mathrm{Tr}(D^{-z}Ag)$. This function is a priori defined only on a
half plane $Re(z)>N$, with $N$ large (see below).

We prove that this new ``zeta" function $z\mapsto
\mathrm{Tr}(D^{-z}Ag)$ has a meromorphic extension to the complex
plane with only simple poles. These poles can be only at $z = d/m,
(d-1)/m , (d-2)/m ,\ldots$, where
$d = \textrm{ord}(A) + \dim(M^g)$ and $m = \textrm{ord}(D)$.
Furthermore, since all the poles are simple for any conjugacy
class $\langle\gamma\rangle\in\langle \Gamma\rangle$, the sum of
residues
$$
Tr_R^{\langle\gamma\rangle}(A) =\sum_{g\in\langle\gamma\rangle}
res_{z=0}Tr(D^{-z}Ag),
$$
defines a trace on $\maA(M)\rtimes \Gamma$ that is independent of
choice of the invariant, positive operator $D$. This is the first
version of our equivariant noncommutative residue.

Given an irreducible representation $\pi$ of the group $\Gamma$, the
asymptotic of the counting function $N_{\pi,D}$ for an order $m$
operator $D$ can now be obtained by applying the standard Tauberian
theorem to the functions $z\mapsto Tr(D^{-z}g)$. One can then prove
that if the group action on $M$ is faithful, then
$N_{\pi,D}(\lambda)$ grows asymptotically as $C \dim(\pi)
\lambda^{n/m}$, where $n=\mathrm{dim}(M)$ and $m=\mathrm{order}(D)$.

Wodzicki shows that the noncommutative residue of an operator $P$ can
be locally be obtained by integrating the $-n$ component
$\sigma_{-n}(P)$ of the asymptotic expansion in local coordinates,
over the cosphere bundle $\operatorname{Tr}_R(P) =
\int_{S^*M}\sigma_{-n}(P)$ with respect to volume measure coming from
the symplectic structure on $T^*M$. This is remarkable because the
local expansion is not diffeomorphism invariant.  A similar formula
for the equivariant case will involve not just the $-n$th term in the
local asymptotic expansion, but also normal derivations over fixed
point submanifolds of higher order (i.e $>-n$) terms from the
expansion (see Proposition \ref{local_form}.)

Also we can now work back on Connes Trace formula that tells us that
the noncommutative residue is an extension of the Dixmier trace to all
pseudodifferential operators. We thus obtain another equivalent
construction for of the Dixmier trace of an invariant operator.
  
We begin with the analysis of the zeta functions $z\mapsto
\operatorname{Tr}(D^{\!-\!z}Ag)$. The proof that these functions are
meromorphic requires a careful use of the stationary phase
principle. We prove, in fact, a slightly stronger version of the
stationary phase principle needed to obtain our local formula for the
noncommutative residues. Finally we use the Tauberian theorem to
obtain the asymptotics of representations in eigenspaces.
        
\subsection*{Acknowledgment} I would like to thank my supervisor Victor Nistor
for all his suggestions and comments and  above all his patience with my work.

\section{Explicit description of Traces}\label{section-EDT}

For a finite group $\gp$ acting on a closed manifold $M$, the
description of the Hochschild homology of the crossed product algebra
$\maA(M) \rtimes \gp$ determines the dimension of the space of traces on
the cross-product algebra $\maA(M)\rtimes\gp$ (see \cite{Dave}).
These traces correspond to conjugacy classes $\langle \gamma \rangle$
of the group $\gp$, for when the elements $\gamma\in \langle \gamma
\rangle$ have a nonempty fixed point set on the cosphere bundle or
$S^*M^{\gamma}\neq\emptyset$. We wish to describe these traces more
explicitly.

By a constant order holomorphic family $A(z) \in \Psi^m(M)$ of
pseudodifferential operators on a compact, smooth manifold $M$, we
mean a family that can be obtained from a holomorphic family of
complete symbols of fixed order and a holomorphic family of
regularizing operators. (The regularizing operators are given the
Fr\'echet topology of $\mathcal{C}^\infty(M \times M)$.) By a
holomorphic family $A(z) \in \Psi^z(M)$, 
we mean that $A(z) = B(z) D^{z} + C(z)$, where $B$ is a holomorphic
family of order zero operators, $D$ is a positive first order
pseudodifferential operator, and $C(z)$ is a holomorphic family of
regularizing operators. See \cite{ALNV, LauterMoroianu1,
LauterMoroianu2, Moroianu} and the references therein for more
details.

Let us fix a positive, elliptic, order one operator $\maD\in\Psi^1(M)$
that is invariant under $\gp$ action. This is to say $\maD g=g\maD$ as
operators on $\smooth(M)$ for all $g\in\gp$.  The complex powers
$\maD^z$ can be defined easily by spectral theorem, but it is much
harder to prove that each $\maD^z$ is a pseudo-differential operator
of order $z$. This was proved by {Seeley} \cite{Seeley}.  Another
proof was given by {Guillemin} in \cite{Guillemin}.

In particular, Seeley result implies that for $Re(z)< -n$, $\maD^z$ is
a trace class, and hence the map
\[z\rightarrow tr(\maD^z)\] 
is holomorphic on the half plane $Re(z)< -n$ of complex numbers.

For $A$ in $\Psi^{\infty}(M)$ and any group element $g$, we define
\[\zeta_{g,A}(z):=\operatorname{Tr}(\maD^{\!-\!z} Ag).\]
This function is a priori defined only when $\maD^{-z}A$ is a trace
class, that is, when $Re(z)>n+order(A)$. We wish to establish the
meromorphic continuation of the function $\zeta_{g,A}(z)$. This would
involve the use of the stationary phase principle. We recall the
following two theorems from \cite{Hormander}.

\begin{theo}\label{rapidly}
For an open domain $\maU\subset \RR^n$, let $u(x)\in
\smooth_c(\maU)$. Let $f\in \smooth (\maU)$ be such that $f'(x)\neq 0
$ in $Supp(u)$. Then for every $ k\in \mathbb{Z}$ we have
\[ 
\left|\int e^{irf(x)}u(x)dx\right| <r^{-k}\] for $r$ large
enough. That is, the integral is a rapidly decreasing function of $r$.
\end{theo}

Let $\maU\subset \RR^n$ be an open domain. Let $f\in \smooth (\maU)$
and $x_0$ be an isolated non-degenerate critical point of $f$ such
that $f'(x_0)= 0\,,\,Im f(x_0)\geq 0$ and $f''(x_0)\neq 0$. Let
\[g(x)=f(x)-f(x_0)-\langle f"(x_0)(x-x_0),x-x_0\rangle\]
and for $j\geq 0$ define the operators,
\begin{equation}\label{the_ljs}
	L_j(u) := (det(f"(x_0)/2\pi i))^{-\!\frac{1}{2}}
	\sum_{p-q=j}\sum_{2p\geq3q}i^{-j}2^{-p}\langle
	f"(x_0)^{-\!1}\maD,\maD\rangle^p(g^qu).
\end{equation}
The operators $L_j$ are order $2j$ constant coefficient operators, for
a given function $f(x)$.

\begin{theo}[Stationary Phase Principle]\label{Stationary}
For an open domain $\maU\subset \RR^n$, let $u(x)\in \smooth_c(\maU),
f\in \smooth (\maU)$, and $r>0$. If the point $x_0\in\maU$ is such
that $f'(x_0)= 0\,,\,Im f(x_0)\geq 0$, and $\operatorname{det}f''(x_0)\neq 0$ is the
unique stationary point of $f$ in $Supp(u)$, then there exist
constants $M_j$ such that, for every $ k\in \mathbb{Z}$, 
\[ \left|\int e^{irf(x)}u(x)dx-\sum_{j<k}e^{irf(x_0)}
M_jr^{-\frac{n}{2}-j}\right|< r^{-k}\] for $r$ large enough.
Furthermore, each $M_j$ is of the form $M_j=L_j(u)(x_0)$, where $L_j$
is a differential operator defined by \eqref{the_ljs}.
\end{theo}

Thus a function $\rho(r)=\int e^{irf(x)}u(x)dx$, where $u$ and $f$
satisfy the above conditions, is asymptotically of order
$-\frac{n}{2}$ in the highest term as $r\rightarrow \infty$. The
following corollary states the same result for a holomorphic family
of such integrals.

\begin{cor}
If $u^z(x)$ is a holomorphic family of functions in $\smooth_c(\maU)$,
and $x_0$ is the only critical point of $f$ in $supp(u^z)$. Then for
$z$ in a compact set, there exist holomorphic functions $M_j(z)$ and a
constant $C>0$, such that \[ \left|\int e^{irf(x)}u^z(x)dx -
\sum_{j<k}e^{irf(x_0)}M_j(z) r^{-\frac{n}{2}-j}\right| < r^{-k}\]
for $r$ large enough.
\end{cor}

\begin{proof}
By the Theorem \ref{Stationary}, we can define
$M_j(z)=L_j(u^z)(x_o)$. Now we only have to check that $M_j(z)$ is
holomorphic. Since $u^z$ is a holomorphic family of functions
$z\mapsto u^z(x_0)$ is a holomorphic function.  By the definition of
holomorphic family, $\partial_ju^z$ is a holomorphic family of
functions, and therefore each $L_j(u^z)$ is, since $L_j$ are constant
coefficient operators. Thus $M_j(z)$ are holomorphic functions as
desired.
\end{proof}

We are now ready to extend $\zeta_{g,A}$ to a meromorphic function.
For any fixed $g \in \gp$, we shall denote by $M^g$ the set of fixed
points of $g$ acting on $M$.

\begin{prop}\label{meromorphic}
Let $k_g=\dim(M^g)$ and $\maD$ be a positive definite elliptic
operator invariant under $\gp$. Then the function $z\mapsto
\operatorname{Tr}(\maD^{\!-\!z} Ag)$ is holomorphic on the half
plane $Re(z)>d,\,\,d=k_g+order(A)$, and has a meromorphic extension
with possible simple poles at $z = d, d-1, d-2, \ldots$.
\end{prop}

\begin{proof}
For the purpose of simplifying the notation, we say that two
functions over an open domain in the complex plane are equivalent
$f\simeq g$ if the difference $f-g$ extends to an entire holomorphic
function.

Given any neighborhood of the diagonal $U$, the operator
$\maD^{-\!z}A$ can be written as a sum of $\maD^{-\!z}A=P(z)+S(z)$,
where the operator $P(z)$ has its kernel supported in $U$ and $S(z)$
is a family of smoothing operators. Since $\operatorname{Tr}(S(z)g)$
is a holomorphic function, the zeta function
$\operatorname{Tr}(\maD^{-\!z}Ag)\simeq
\operatorname{Tr}(P(z)g)$. Thus we could assume without any loss of
generality that the family $\maD^{-\!z}A$ is supported sufficiently
close to the diagonal.

For $Re(z)> n+order (A)$, the operators $\maD^{-z}A$ are trace
class, and hence the map $z\mapsto \operatorname{Tr}(\maD^{\!-\!z}
Ag)$ is holomorphic on the half plane $Re(z)>n+order(A)$. We wish to
show that it has a meromorphic continuation.

Let us fix a $\gp$ invariant metric on $M$. For every point $x\in M$,
choose a normal coordinate chart $\maU_x$. We can assume that the
closures $\overline{\maU_x}$ and $\overline{g\maU_x}$ are either equal
or disjoint. Thus, when $x\in M^g$, we assume that $U_x$ is $g$
invariant. We choose a finite subcover $\maU_i=\maU_{x_i}$ for
finitely many $x_i\in M$. We choose $\psi_i\prec \maU_i$ a partition
of unity subordinate to this cover.

For any trace class pseudo-differential operator $P$,
\begin{align*}
  \operatorname{Tr}(Pg)&=\sum_{j,k}
  \operatorname{Tr}(\psi_kP\psi_jg)\\ &=\sum_{j,k}
  \operatorname{Tr}(\psi_j ^{\frac{1}{2}}g\psi_kP\psi_j
  ^{\frac{1}{2}}).
\end{align*}
If the kernel of $P$ is supported in the open set $\Omega\supset
\Delta$ such that $\overline{\Omega}\cap\overline{g\maU_i\times
\maU_i}=\emptyset$ whenever $x_i$ is not a fixed point, then the sum
above reduces to sum about the fixed point manifold namely,
\[\operatorname{Tr}(Pg)=\sum_{\maU_j\cap g.\maU_j\neq\emptyset} 
\operatorname{Tr}(\psi_j ^{\frac{1}{2}}gP\psi_j ^{\frac{1}{2}}).\] The
trace of a smoothing family is an entire function, and hence we could
subtract a suitable family of smoothing operators from $\maD^{-z}A$,
and we assume that its supported in a small enough neighborhood of the
diagonal as above. We can thus restrict ourselves to coordinate
neighborhoods of the fixed points.

Similarly, we only need to consider operators of the form
$\psi^{\frac{1}{2}}\maD^{\!-z}A\psi^{\frac{1}{2}}$ when we work on
local coordinates with $\psi_j$ as partition of unity subordinate to
a cover as above. Then $\maD^{-z}A$ restricted to the functions
supported in such a neighborhood $\maU$ could be given in terms of
quantization of a symbol $p^z(x,\xi)\in S^{-z+a}(\maU)$ as
\begin{eqnarray}\label{quantize}
  \dot{q}(p^z)u(x)=\int e^{i\langle x-y,\xi\rangle}p^z(x,\xi)u(y)
  dy d\xi=\psi^{\frac{1}{2}}\maD^{_z}A\psi^{\frac{1}{2}}u(x).
\end{eqnarray}
 Thus by $\gp$ invariance of the metric we get
\begin{align*}
  \psi^{\frac{1}{2}}\maD^{-z}A\psi^{\frac{1}{2}}g(u(x))&
  =\int e^{i\langle x-y,\xi\rangle}p^z(x,\xi)g\cdot u(y) dy d\xi\\
  &=\int e^{i\langle x-gy,\xi\rangle}p^z(x,\xi) u(y) dy d\xi.
\end{align*}
By abuse of notation, we shall continue to denote the operator of the form
$\psi^{\frac{1}{2}}\maD^{\!-\!z}A\psi^{\frac{1}{2}}$ as
$\maD^{\!-\!z}A$.

The distributional kernel of $D^{-z}Ag$ restricted to $\maU\times
\maU$ is therefore given by
\[K(x,y)=\int e^{i\langle x-gy,\xi\rangle}p^z(x,\xi) d\xi.\] Thus when the
operator $\maD^{-z}Ag$ is trace class (for sufficiently large values
of $z$), the kernel $K(x,y)$ is integrable and the trace is given by
\begin{eqnarray}
\label{integral}\zeta_{g,A,\maU}(z)=\int K(x,x)dx=\int e^{i\langle
  x-gx,\xi\rangle}p^z(x,\xi) d\xi dx.
\end{eqnarray}
The integral on the right hand side is absolutely convergent for
$Re(z)>n+a$, and we would like to establish the continuation of this
integral to a meromorphic function on the complex plane.

To apply Stationary phase to the Integral (\ref{integral}), we
notice first that in $\gp$ invariant normal coordinates, $g$ becomes
an orthogonal linear map with $g^*=g^{-\!1}$. Then the derivative of
the phase function $f(x,\eta)= \langle x-gx,\eta\rangle$ would be
\begin{eqnarray}
\label{Df}
Df_{(x,\eta)}(\alpha_1,\beta_1)=\langle x-gx,\beta_1\rangle+\langle
\alpha_1-g\alpha_1,\eta\rangle 
\end{eqnarray}
 that implies that the stationary
points $(x.\eta)$ of the phase function correspond to the points $(x,\eta)$ such
that $x=gx$ and
$\eta=g\eta$. Furthermore,
\begin{eqnarray}\label{D2f}
f''(x,\eta):=D^2f_{(x,\eta)}(\alpha_1,\beta_1)(\alpha_2,\beta_2)=\langle
\alpha_2-g\alpha_2,\beta_1\rangle+\langle
\alpha_1-g\alpha_1,\beta_2\rangle .
\end{eqnarray}
 We choose a special normal
coordinate chart. Since we are only interested in the neighborhoods
of $x\in M^g$, the tangent space here decomposes as $T_xM=T_xM^g\oplus
(1-g)T_xM$. We choose an orthonormal basis for $T_xM^g$ and $(1-g)T_xM$
to determine a normal coordinates at x. Set $\maU^g=exp(T_xM^g)$ and
$(1-g)\maU=exp((1-g)T_xM)$. That is,
\[\maU\ni x=(x_1,x_2)\in \maU^g\times (1-g)\maU.\]
After trivializing the the cotangent bundle $T^*\maU=\maU\times \RR^n$ one can also rewrite
\[\RR^n=\RR^{n^g}\times (1-g)\RR^n,\]
and henceforth $\xi=(\xi_1,\xi_2)$. Set $\pi_x(x)=x_1$ and
$\pi_{\xi}(\xi)=\xi_1$. And for each $(x',\xi')\in T^*\maU^g$, set $
F_{(x',\xi')}:=\pi_x^{-1}(x')\times \pi_{\xi}^{-1}(\xi')\subset
T^*\maU$.

Now the symbol $p^z(x,\xi)\in S^{-z+a}(\maU)$ has an asymptotic
expansion of the form
\[p^z(x,\xi)\simeq p^z_{-z+a}(x,\xi)+p^z_{-z+a-1}(x,\xi)+\ldots
+p^z_{-z+a-j}(x,\xi)\ldots, \] where each symbol $p^z_l(x,\xi)$ is
homogeneous of degree $l$ for $|\xi|\geq 1$. By the holomorphicity of
the family $\maD^{-\!z}A$, each $p^z_l(x,\xi)$ is also holomorphic
in $\smooth(T^*U/\{0\})$. Further this expansion determines the
operator $\maD^{-z}A$ up to a smoothing operator. If the $Re(z)$ is
in the interval $[J-1,J)$, then the operator
\[Q(z):=\dot{q}\left(p^z(x,\xi)-\sum_{j\leq J+n}p^z_{-z+a-j}(x,\xi)\right),\]
 is trace class and hence the function
$\operatorname{Tr}(Q(z)g)$ is entire. Thus
$$\zeta_{g,A,\maU}(z)\simeq
\sum_{Re(l)\geq-n}\operatorname{Tr}(\maD^{-\!z}\dot{q}(p^z_l)g).$$

Now it suffices to check that $z\mapsto
\operatorname{Tr}(\maD^{-\!z}\dot{q}(p^z_l)g)$ has a meromorphic
continuation for the operators defined by the each term in
right-hand side, that is, for $l=-z+a,-z+a-1,-z+a-2\ldots$ so on
till $Re(l)\geq -n$. (Because as we have just shown only finitely many of these
terms contribute tot he trace.) For this one must show that the integral in
(\ref{integral}) for $\dot{q}(p^z_l)$ has a meromorphic
continuation. Since this integral on any compact subset in $\xi$
would yield an entire function, we could restrict the domain where
$p^z_l$ is homogeneous, which is $|\xi|>1$.

Let $r=|\xi|$. Let $dS^n(r)$ be the volume form on a sphere of
radius $r$. Then
\begin{align*}
\zeta_{g,\dot{q}(p^z_l),\maU}(z)&=\int_{|\xi|>1} e^{i\langle
    x_2-gx_2,\xi\rangle} p^z_l(x_1+x_2,\xi) d\xi dx_2\,dx_1\\
  &+\text{entire function}\\
  &\simeq \int_1^{\infty}\int_{|\xi|=r} e^{i\langle
    x_2-gx_2,\xi_2\rangle}p^z_l(x_1+x_2,\xi)dS^n(r)dxdr.
\end{align*}

Let $dS^n=dS^n(1)=r^{1-n}dS(r)$ be volume form on the unit sphere in
$\RR^n$. We would also use the notation $dS^n(\xi)$ to emphasize the
variable used in integration. From the homogeneity of $P^z_l$ we
obtain
\begin{multline*}
  \zeta_{g,\dot{q}(p^z_l),\maU}(z) \simeq
  \int_1^{\infty}r^{n+l-1}\int_{|\xi|=1} e^{ir\langle
    x_2-gx_2,\xi_2\rangle} p^z_l(x_1+x_2,\xi)dS^n(\xi)dxdr\\
  \simeq\int_1^{\infty}r^{n+l-1}\int_{|\xi_1|< |\xi_2|}\ e^{ir\langle
    x_2-gx_2,\xi_2\rangle}p^z_l(x_1+x_2,\xi)
  dS^n(\xi) dx_2dx_1dr\\
  \qquad +\int_1^{\infty}r^{n+l-1}\int_{|\xi_1|\geq|\xi_2|}\ e^{ir\langle
    x_2-gx_2,\xi_2\rangle}p^z_l(x_1+x_2,\xi)
  dS^n(\xi) dx_2dx_1dr .
\end{multline*}
Set 
\[\zeta_{g,\dot{q}(p^z_l),\maU}(z)\simeq I_1(p^z_l,z)+I_2(p^z_l,z)\]
  
The first integral $I_1(p^z_l,z)$ above yields an entire function.
On the set $|\xi_1|<|\xi_2|$, we have $|\xi_2|>0$ and thus  by \ref{Df} the phase
function of the integral
\[I_1(p^z_l,z)=\int_{|\xi_1|< |\xi_2|}\ e^{ir\langle x_2-gx_2,\xi_2\rangle}
p^z_l(x_1+x_2,\xi) dS^n(\xi) dx_2dx_1\] will not have any critical
points. Hence by Theorem \ref{rapidly}, it must be rapidly
decreasing asymptotically as $r\rightarrow \infty$. After
integrating in $r$, we have $I_1(p^z_l,z)$ converges for all $z$ and
gives an entire function.

However, the phase function $f$ of the integral
\begin{eqnarray}\label{main_integral}
I_2(p^z_l,z)=\int_{|\xi_1|\geq|\xi_2|}\ e^{ir\langle
x_2-gx_2,\xi_2\rangle}p^z_l(x_1+x_2,\xi) dS^n(\xi) dx_2dx_1
\end{eqnarray}
has a stationary phase point whenever $x_2=0$ and $\xi_2=0$, that is
indeed the fixed point manifold $S^*\maU^g$. For $\xi=(\xi_1,\xi_2)$
on the unit sphere with $|\xi_1|\geq|\xi_2|$, let
$\eta_1=\frac{\xi_1}{|\xi_1|}\in S^{k_g}(\xi_1)$ and
$\eta_2=\frac{\xi_2}{|\xi_1|}$. Then rewriting (\ref{main_integral})
above we have
\[\int\left(\int_{|\eta_2|\leq 1}\ e^{ir\langle x_2-gx_2,\eta_2\rangle}
  p^z_l(x_1+x_2,\eta_1+\eta_2) d\eta_2 dx_2\right)dS^k(\eta_1)dx_1dr.\] But on
each normal fiber $F_{(x_1,\eta_1)}$ of a fixed point
$(x_1,\eta_1)\in S^*M^g$, there is only one stationary point of $f$
namely $(x_1,\eta_1,0,0)$ and the Hessian $D^2f\neq 0$. Then by the
stationary phase method \ref{Stationary}, for every integer $i\geq 0$ there exist
$\overline{M}_j(z)(x_1,\eta_1)$ $0\leq j\leq i$, smooth in $x,\xi\in S^*\maU$ and
holomorphic in $z$, and a constant $c_i>0$ such that
\begin{multline*}
  \left|\int_{|\eta_1|\geq|\eta_2|}\ e^{ir\langle
    x_2-gx_2,\eta_2\rangle}p^z_l(x_1\!+\!x_2,\eta_1\!+\!\eta_2) d\eta_2
  dx_2\!-\!\sum_j \overline{M}_j(z)(x_1,\eta_1) r^{-n+k_g-j}\right|\\
  <C_ir^{-i}\qquad \text{as}\,r~\rightarrow \infty.
\end{multline*}
Set
\begin{eqnarray}\label{asymp_estimate}
M_j(z)=\int \overline{M}_j(z)(x_1,\eta_1)dx_1 d\eta_1.
\end{eqnarray}

Therefore, up to a holomorphic function,
\begin{eqnarray*}
I_2(p^z_l,z)&\simeq &\int_1^{\infty} \sum_j
M_j(z)r^{k_g+l-1-j}\,dr\\&=&\sum_j \frac{ M_j(z)}{k_g+l-j}=\sum_j
\frac{ M_j(z)}{k_g+a-z-j-i}.
\end{eqnarray*}
The last equality is due to the fact that $l$ is the degree of
homogeneity, and takes value in the set $a-z-i\,,\, (i\in \ZZ)$.
Thus $I_2$ and hence $\zeta_{g,A}$ would be meromorphic with at most
simple poles as specified.
\end{proof}

\begin{remark}
In the proof of Lemma \ref{meromorphic} above, the operators
$\maD^{-z}A$ could be replaced by any holomorphic family $R(z)$ with order
 $R(z)=-z+a$.
\end{remark}

\begin{remark}
For any order $m$ operator positive $D$, the operator
$\maD=D^{\frac{1}{m}}$ is of order one and hence for the operator $D$,
its zeta function $\zeta_{g,A}(z)=\operatorname{Tr}(D^{-\!z}Ag) $
would have a meromorphic extension to the complex numbers with
possible poles at $\frac{d}{m},\frac{d-1}{m},\ldots.$
\end{remark}

\begin{remark} Nonsimple poles can occur in
the zeta function of operators on spaces with
singularities \cite{Moroianu}.
\end{remark}

Although the function $\zeta_{g,A}$ does depend on the choice of an
invariant operator $\maD$, the residue of $\zeta_{A,g}$ at $z=0$ does
not. Indeed if $\maD_1$ is another such operator invariant under
$\gp$, then let $R(z)=\maD_1^{-z}\maD^z.$ The family $R(z)\in
\Psi^0(M)$ and $R(0)=I$. Thus we define
\begin{equation}\label{def.residue}
\tau_g(A) :=
res_{z=0}\operatorname{Tr}(\maD^{-z}Ag)=
res_{z=0}\operatorname{Tr}(R(z)\maD^{-z}Ag)
=res_{z=0}\operatorname{Tr}(\maD_1^{-z}Ag).
\end{equation}
 
It is also immediate from the definition that the residue $\tau_g$
defined in Equation \eqref{def.residue} descends to a function on
$\maA(M)=\pdoM$.

In particular, we obtain the proof of Proposition \ref{meromorphic}
that the residue vanishes if $g$ has only isolated fixed points.

\begin{cor}
If $g$ has no fixed points on the co-sphere bundle $S^*M$, then
$\tau_g=0$.
\end{cor}

This can explain in a certain way the result of Atiyah and Bott
in \cite{AB1}.

\begin{cor} Let $\maB(M)=\maA(M)\rtimes \gp$.
For any conjugacy class $\langle \gamma \rangle \in \langle \gp
\rangle$, the map $\operatorname{Tr}_R^{\langle \gamma
\rangle}:\maB(M)\mapsto \CC$ given by
\[
\operatorname{Tr}_R^{\langle \gamma \rangle}
(\sum_{g\in \gp}A_gg):=\sum_{g\in\langle 
\gamma \rangle} \tau_g(A_g)
\] 
is a trace on $\maB(M)$.
\end{cor}

\begin{proof}
By linearity it suffices to check that each $\tau_{\gamma}$ is a trace
on the generators. If $g,h\in \gp$ then $gh$ and $hg$ belong to same
conjugacy class and thus
\[\tau_{\ip{gh}}[Ag,Bh] =res_{z=0}\operatorname{Tr}(D^{\!-\!z}[Ag,Bh]).\]
For any operators $A,B,C$, such that $ABC$, $ACB$, and $BAC$ are trace
class and $\operatorname{Tr}(ACB)=\operatorname{Tr}(BAC)$ the
following obviously holds:
\[\operatorname{Tr}(A[B,C])=\operatorname{Tr}([A,B]C).\]
Thus for $Re(z)$ large enough,

\[\operatorname{Tr}(\maD^{-z}[Ag,Bh]) =
\operatorname{Tr}([\maD^{-z},Ag]Bh) =
\operatorname{Tr}([\maD^{-z},A]gBh). \] Here the  equalities holds
because $\maD$ commutes with $g$. But then $[\maD^{-z},A]=zR(z)$, for
some holomorphic family $R(z)$.  Since $\operatorname{Tr}(R(z))$ can
at most have a simple pole at $z=0$, $res_{z=0}zR(z)=0$.
\end{proof}

\begin{remark} If the 
elements in $\ip{\gamma}$ have no fixed points on the cosphere bundle
$S^*M$ then $\operatorname{Tr}_R^{\ip{\gamma}}=0$ as the zeta function
of these elements have a holomorphic continuation to the complex
plane. But here we remark that by an argument as above the sum of
regular values $Ag\rightarrow \zeta_{A,g}(0)$ defines a trace on
$\Psi^{\infty}(M)\rtimes \gp$. This trace will still be non-trivial if
$g$ has isolated fixed points on $M$, but will not descend to be a
trace on $\maA(M)\rtimes \gp$.
\end{remark}

For simplicity, let us assume that $S^*M^g$ is connected or
empty for the rest of our discussion.  Then the Hochschild homology
calculation in \cite{Dave} show that all the traces of the algebra
$\maB(M)$ are linear combinations of the traces
$\operatorname{Tr}_R^{\ip{\gamma}}$.

\begin{lem}\label{nontrivial}
Let $m=-k_g$, and let $A$ be of order $m$ with principal symbol
$a_m=\sigma(A)$. Then there exists a constant $C>0$ such that
$\tau_{g}(A)=C\int_{S^*M^g} \sigma(A)dvol$.

In particular, if $S^*M^g\neq\emptyset$ then
$\tau_{}(\maD^{-k}g)\neq 0$ for
  a positive invariant order one operator $\maD$.
\end{lem}

\begin{proof}
Since the principal symbol of an operator is defined independent of
coordinates, we can work with normal local coordinates. By Proposition
\ref{meromorphic}, the first possible pole of $\zeta_{g,A}(z)$ would
be at $z=0$. In each normal local coordinate, the value of the residue
is given by $M_0(0)$, which can be obtained by the Equation
(\ref{asymp_estimate}) \[M_0(z)=\int
M_0(z)(x_1,\eta_1)dS^k(\eta_1)dx_1.\] By the stationary phase
principal \ref{Stationary}, we have
$$M_0(0)(x_1,\xi_1)=L_0(p^0_{-k_g}(x_1,\xi_1))=C\sigma(A)(x_1,\xi_1),$$
where $C=\textrm{det}(D^2f_(x_1,\xi_1)/2\pi)^{-\frac{1}{2}}$ is a
constant independent of the fixed point $(x_1,\xi_1)$. Since the
principal symbol $\sigma(A)$ is invariantly defined under change of
coordinates, the local formula is independent of coordinates and adds
up to give the result.
\end{proof}

Related calculations (but not in the equivariant setting)
can be found in several other papers, including
\cite{LauterMoroianu1, LauterMoroianu2, Lescure, 
MoroianuNistor, Schrohe}.

\section{Local Expression for  the traces}\label{section-LET}

We now establish a local formula for the traces
$\operatorname{Tr}_R^{\ip{\gamma}}$ for all operators. This in
particular would require a generalization of Lemma \ref{nontrivial}
and a more elaborate use of the stationary phase principal. For any
operator $A$ of order m, the calculation in Proposition
\ref{meromorphic} yields a local formula for the residue of
$\zeta_{g,A}$ at $z=0$. To start, fix a $\gp$ invariant metric on $M$. By
compactness of $M$ let $r$ be the minimum value of injectivity radius on $M$. Let
$\maU$ be a tubular neighborhood of $M^g$ such that $dist(y,gy)<\frac{r}{2}$ for
any $y\in \maU$. Denote by $j(x)$ the tangent at $x$  to the unique geodesic
$\gamma$ such that $\gamma(0)=x$  and $\gamma(1)=gx$. Then we define a  function
on $T^*\maU$ by
\[f(x,\xi):=\ip{\xi,j(x)}.\]
 In fact locally in normal coordinates $f(x,\xi)=\ip{x-gx,\xi}$. Thus
 Equation \ref{Df} implies that $M^g$ is the critical manifold of $f$ and by
 Equation \ref{D2f} $f$ is Morse Bott, the Hassian $f''=:D^2f$ is a non-degenerate
 bilinear form on the normal bundle
 $\operatorname{N}(T^*\maU:T^*M^g)=(1-g)T^*\maU_{|T^*M^g}$. Therefore the  inverse
 $D^2f^{\!-\!1}$ is a  nondegenerate bilinear form on the conormal bundle
 $\operatorname{N}^*(T^*\maU:T^*M^g)$.

 Let
$\maU=(x_1,x_2)$ be normal local coordinates
at a fixed point $x\in M^g$ and choose
$\eta_1,\eta_2$ to give a trivialization over $\maU$ of the cotangent
bundle $T^*M^g$ and its orthogonal $(1-g)T^*M$ in the cotangent
bundle.  Then by the Equation (\ref{asymp_estimate}),
\[res_{z=0}\operatorname{Tr}(\psi^{\frac{1}{2}}
\maD^{-\!z}Ag\psi^{\frac{1}{2}})=\sum_{j=0}^{m+k_g}\int
\overline{M}_j(x_1,\eta_1) dS^{k_g}(\eta_1)dx_1.\] Let $D$ denote the
gradient operator on the normal bundle $(1-g)T^*M$, that is,
$D^{\perp}=i\partial_{x_2}+i\partial_{\eta_2}$. Again by stationary phase
principle, each $M_j$ can be calculated as
\begin{align*}
\overline{M}_j(0)(x_1,\xi_1) &= L_j(p^0_{\!-\!k_g+\!j}(x_1,\eta_1))\\
&\\ &=C\langle f''^{\!-\!1}_{(x_1,\eta_1)} D^{\perp},D^{\perp}\rangle^j
\left(p^0_{\!-\!k_g\!+\!j}(A)(x_1,\eta_1)\right),
\end{align*}
where the constant $C=\left(\textrm{det}\frac{f''_{(x_1,\eta_1)}}{2\pi
i}\right)^{-\frac{1}{2}}$.

Now by using the invariance of the residue under the choice of of
operator defining the $\zeta_{g,\maA}$ one has:

\begin{prop} \label{local_form}
Let $\maU(x_1,x_2)$ be normal coordinates near a point $x\in M^g$ and
let $A\in\Psi^m(M)$ be an order $m$ operator given on $\maU$ by the
quantization of the symbol $a(x,\eta)\sim
a_m(x,\eta)+a_{m-1}(x,\eta)+\ldots$ Let $dS^{k_g}(\eta)$ denote the
volume form on $k_g-1$ dimensional sphere in $S^*\maU^g\simeq
\maU\times S^{k_g-1}$.

Let
\begin{align*}
W_{M^g}(A):&=\left(\sum_{j=0}^m\int_{S^k(\eta_1)} C\langle
  f''^{\!-\!1}_{(x_1,\eta_1)}D^{\perp},D^{\perp}\rangle^ja_{-k_g+j}(x,\eta)
  dS^{k_g}(\eta_1)\right)|dx_1|
\end{align*}
Then $W_{M^g}(A)$ is invariant under $g$-equivariant  change
of normal coordinates near $M^g$.
\end{prop}

\begin{proof}
Let $\phi:M\rightarrow M$ be an $\gp$ equivariant isometry. Then
$\phi(M^g)=M^g$.  Let $A$ be supported in a normal neighborhood $\maU$
as above and let $\phi(\maU)=\maV$.  Then from the discussion above
(applied to $\phi_*A$) it follows that:
\begin{align*}
W_{M^g}(A)&=\operatorname{res}_{z=0} \operatorname{Tr}(\maD^{-z}\phi A
\phi^{\!-\!1}g)\\ &=\operatorname{res}_{z=0}
\operatorname{Tr}(\phi^{\!-\!1}\maD^{-z}\phi A g)
\end{align*}
Since $\phi$ is an isometry,
\[
\phi^{\!-\!1}D^{\!-\!z}\phi=\phi^*D^{\!-\!z}\phi=\maD_1^{\!-\!z},
\]
Where $\maD_1$ is another order $1$ positive elliptic operator. Since
the residue is independent of choice of such an operator, we have

\[W_{M^g}(\phi_*A)=\operatorname{res}_{z=0}\operatorname{Tr}(\maD^{\!-\!z}_1
Ag)=W_{M^g}(A).\]
\end{proof}

\section{Asymptotics for representations in
eigenspaces}\label{section-ARE}

Let $\pi\in R(\gp)$ be a representation of the group $\gp$. We can
define a trace on the cross-product algebra $\maB(M)$ corresponding
to $\pi$ as
 \[ \tau_{\pi}\left(\sum_{g\in\gp}A_gg\right):=\frac{1}{|\gp|}\sum_{g\in
\gp}\overline{\chi_{\pi}(g)}\tau_g(A_gg)\]. 

Since $\chi_{\pi}$ is class function, this trace is a linear
combinations of the traces $\operatorname{Tr}_R$ defined before. Lets
denote by $\epsilon:=\frac{1}{|\gp|}\sum_{g\in\gp} g$ the idempotent
in $\CC[\gp]$. Then one can likewise obtain traces on the invariant
algebra $\maA(M)^{\gp}\,$ by $A\rightarrow
\tau_{\pi}(A\epsilon)$. Such a trace is the residue of a $\zeta$
function namely the one obtained from the meromorphic extension using
Proposition \ref{meromorphic} of
\[\zeta_{\pi,A}(s)=\frac{1}{|\gp|}\sum_{g\in
\gp}\overline{\chi_{\pi}(g)}\operatorname{Tr}(\maD^{-\!s}Ag)\qquad
Re(s)>n.\] Here $\maD$ is $\gp$ invariant positive order one operator.

We say that action of $\gp$ on $M$ is faithful or effective if for any
$g\in\gp$ the fixed point manifold $M^g=M$ then $g$ must be the
identity element.

\begin{lem}\label{effective}
Let $n=\textrm{dim}(M)$. If the action of $\gp$ on $M$ is effective,
then $\zeta_{\pi,Id}$ has a pole at $s=n$.
\end{lem}

\begin{proof}
Let $k_g=\dim{M^g}$ for any group element $g$. If the action of the
group $\gp$ is effective then $k_g<n$ whenever $g$ is not the identity
element $e$. Since $\zeta_{e,Id}(s)$ does have a pole at $s=n$ and all
the $\zeta_{g,Id}$ do not, $\zeta_{\pi,Id}$ is holomorphic in the half
plane $Re(s)>n$ and must have a pole at $s=n$.
\end{proof}

For the rest of section, we assume that the $\gp$ action is
effective. Let $\{\lambda_i\}$ be the set of eigenvalues of $\maD$.
Then each eigenspace $V_{\lambda_i}=ker(\maD-\lambda_i)$ is
invariant under $\gp$, and so acquires a representation
$\pi_{V_{\lambda_i}}$ of $\gp$. We count the multiplicity of 
representation $\pi$ asymptotically in the eigenspaces
$V_{\lambda_i}$. Let \label{counter}
\[N_{\pi,\maD}(\lambda) := \sum_{\lambda_i\leq\lambda}
\langle V_{\lambda_i}(g),\pi\rangle =
\sum_{\lambda_i\leq\lambda}\frac{1}{|\gp|}\sum_{g\in\gp}
\overline{\chi_{\pi}(g)}\chi_{V_{\lambda_i}}.\]

Let $m_i=\langle V_{\lambda_i},\pi\rangle$ denote the multiplicity of
occurrence of $\pi$ in the $V_{\lambda_i}$ if $\pi$ is
irreducible. For the purpose of establishing the asymptotics of the
function $N_{}$, we recall the following Tauberian theorem.

\begin{theo}[Tauberain Theorem]\label{tauberian}
Let $v(x)$ be a non-decreasing function with $v(x)=0$ for all $x <1$.
Suppose the function \[f(z)=\int_1^{\infty}x^{-z}dv(x)\] exists for
$Re(z)>1$ and is analytic, and suppose that there is a constant $A$
such that the analytic function \[F(z)=f(z)-\frac{A}{z-1}\] extends
continuously to the closed half plane $Re(z)\geq 1$. Then
asymptotically, \[ \frac{v(x)}{x}\rightarrow A \text{, as }
x\rightarrow \infty.\]
\end{theo}

\begin{cor}[Asymptotic estimates]\label{rep_asymptotics}
Let the group action of $\gp$ on $M$ be faithful. Then for an
irreducible representation $\pi$ and an $\gp$ invariant operator
$\maD$, the multiplicity counting function, $N_{\pi,\maD}(\lambda)\sim
C \dim(\pi) \lambda^{n/m}$, where $C>0$ and
$m=\operatorname{ord}(\maD)$.
\end{cor}

\begin{proof}
First we notice that for sufficiently large values of the complex
variable $s$, the zeta functions have a nice expression of the form
\begin{align}\label{trace_exp}
\zeta_{g,Id}(s)&=\operatorname{Tr}(\maD^{-s}g) = 
\sum_i\operatorname{Tr}(\maD^{-s}g_{|_{V_{\lambda_i}}})\\
& =\sum_i \chi_{V_{\lambda_i}}(g)\lambda_i^{-s},\qquad Re(s)>k_g.
\end{align}
Then adding them all up, we get
\begin{align}
\zeta_{\pi,Id}(s)&=\frac{1}{|\gp|}\sum_{\gp}\overline{\chi_{\pi}(g)}
\sum_i \chi_{V_{\lambda_i}}(g)\lambda_i^{-s}\qquad Re(s)>k_g\\
&=\sum_i\lambda^{-s}_im_i\\
&=\int_0^{\infty}\lambda^{-s}dN_{\pi,\maD}(\lambda).
\end{align}

Since by Lemma \ref{effective} the function $\zeta_{\pi,Id}(s)$ is
holomorphic on the half plane $Re(s)>n$ and has a non-zero residue $K$
at $s=n$, it follows that $\zeta_{\pi,Id}(s)-\frac{K}{s-n}$ has a
continuation to the closed half line $Re(s)\geq n$. Let
$C=\operatorname{Tr}_R(\maD^{-n})$, then we note that $K$ is of the
form $K =C \frac{\chi_{\pi}(e)}{|\gp|}$. The estimate follows by direct application
of Theorem \ref{tauberian}.
\end{proof}

One could also employ Proposition \ref{meromorphic} again to obtain
relative asymptotic estimates for a pair of irreducible
representations $\pi_j,\,j=1,2$ of same dimension, and compare their
occurrence of in eigenspaces of our invariant operator $\maD$
assuming that for all but a finitely many eigenvalues $\lambda_i$,
$$\langle V_{\lambda_i},\pi_1\rangle\geq \langle
V_{\lambda_i},\pi_1\rangle.$$

\begin{cor}
Let $\pi_j$ be irreducible representations as above, and let
\[N_{\pi_1,\pi_2,\maD}(\lambda):=\sum_{\lambda_i\leq \lambda}\langle
V_{\lambda_i},\pi_1-\pi_2\rangle.\] Let
$k=max_{g\in\gp}\{\textrm{dim}(M^g)|\chi_{\pi_1}(g)\neq
\chi_{\pi_2}(g)\}$. Then asymptotically
$$N_{\pi_1,\pi_2,\maD}(\lambda)\simeq C\lambda^k.$$
\end{cor}

\begin{remark} 
One application for Weyl's theorem is to establish the convergences of
numerical schemes for solving boundary value problems using spectral
decomposition. It can be hoped that the equivariant asymptotic
formulas could be applied to establish rates of convergences for
numerical solution to the boundary value problems which involve some
symmetries but over a domain with not necessarily smooth boundary.
\end{remark}

Fix a $\gp$-representation $\pi$. Now for an $\gp$-invariant compact
operator $\maD$ on $L^2(M)$, let $\{\mu_i\}$ be the decreasing
sequence of eigenvalues counted with multiplicity of $\pi$.  Then each
eigenvalue $\lambda$ is counted the number of times $\pi$ occurs as a
subrepresentation in $\textrm{ker}(\maD-\lambda I)$. We shall denote
this $\pi$-multiplicity of an eigenvalue $\lambda_j$ by $m(\pi)_j$.

\begin{definition} 
We say that a positive invariant operator $\maD$ is $\pi$-measurable
if the following limit exists and is finite
\[\operatorname{Tr}_{\pi}^+(\maD) := 
\lim_{N\rightarrow \infty}\frac{\sum_{i\leq
N}\mu_i}{\textrm{Log}\,N}=\lim_{N\rightarrow \infty}\frac{\sum_{\sum
m(\pi)_j\leq N} m(\pi)_j\lambda_j}{\textrm{Log}\,N}\]
\end{definition}

In particular when $\pi$ is the trivial representation for the trivial
group then $Tr^+_{\pi}=Tr^+$ is the Dixmier trace functional which is
defined on the measurable elements of the Dixmier ideal
$\mathcal{L}^{1,\infty}$. 
 \begin{prop} With the above notation let
$\maD$ be an positive invariant pseudodifferential operator operator
of order $-n$ where $n=\textrm{dim}(M)$ then $\maD$ is $\pi$ measurable and,
\[Tr^+_{\pi}(\maD)=\operatorname{Tr_R}(\maD)=Tr^+(\maD).\] Here of
course the last equality is the Trace formula of
Connes\cite{Connes88}.  \end{prop}

\begin{proof}
In view of Connes' theorem it suffices to show
$\operatorname{Tr}_{\pi}^+(\maD)=\operatorname{Tr}_R(\maD)$. We need
the following lemma

\begin{lem}
Let $\{a_k\}$ be a decreasing sequence of positive numbers such that
$a_k\!\searrow \!0$ and $\sum a_k^s<\infty$ for $s>1$ and $\sum
a^s_k\rightarrow\frac{1}{s-1}$ as $s\!\searrow \!1$. Then
$\sum_{k=1}^N a_k\sim\textrm{Log}\,N$.
\end{lem}

We refer the reader to \cite{Gracia-Bondia} for a proof of the above
lemma.
 
As before we now define 
\begin{align*}
f(s):&=\frac{1}{|\gp|}\sum_{\gp}\overline{\chi_{\pi}(g)}
\operatorname{Tr}(\maD^sg)\\
&=\frac{1}{|\gp|}\sum_{\gp}\overline{\chi_{\pi}(g)} \sum_j
\chi_{V_{\lambda_j}}(g)\lambda_j^s\\
&=\sum_j\lambda^s_jm(\pi)_j=\sum_i\mu_i
\end{align*}

By Proposition \ref{meromorphic} $f(s)$ is analytic in the half plane
$\textrm{Re}(s)>1$ and has a simple pole at $s=1$. Then the
proposition follows from the lemma above.

\end{proof}

\end{document}